\documentclass[10pt,a4paper]{amsart}
\usepackage[utf8]{inputenc}
\usepackage{amsmath}
\usepackage{amsfonts}
\usepackage{amssymb}
\usepackage{mathtools}
\usepackage{paralist}
\usepackage{amsthm}
\usepackage{amscd}
\usepackage{tikz}
\usepackage{graphicx}

\newtheorem{theorem}{Theorem}[section]

\newtheorem{lemma}[theorem]{Lemma}
\newtheorem{proposition}[theorem]{Proposition}

\theoremstyle{definition}

\providecommand{\N}{\mathbb{N}}
\providecommand{\R}{\mathbb{R}}
\providecommand{\Z}{\mathbb{Z}}

\author{Robert Schippa}

\makeatletter
\@namedef{subjclassname@2020}{%
  \textup{2020} Mathematics Subject Classification}
\makeatother

\calclayout
\begin{document}
\title[Resolvent estimates for elliptic operators in three dimensions]{Improved resolvent estimates for con\-stant-\-coef\-ficient elliptic operators in three dimensions}

\keywords{resolvent estimates, Sobolev embedding, Fourier restriction, Limiting Absorption Principle}
\subjclass[2020]{Primary: 35J08, 35J30, Secondary: 46E35}

\begin{abstract}
We prove new $L^p$-$L^q$-estimates for solutions to elliptic differential operators with constant coefficients in $\R^3$. We use the estimates for the decay of the Fourier transform of particular surfaces in $\R^3$ with vanishing Gaussian curvature due to Erd\H{o}s--Salmhofer to derive new Fourier restriction--extension estimates. These allow for constructing distributional solutions in $L^q(\R^3)$ for $L^p$-data via limiting absorption by well-known means. 
\end{abstract}

\maketitle
\section{Introduction}

The purpose of this note is to show new $L^p$-$L^q$-estimates for solutions to elliptic differential equations in $\R^3$. Let
\begin{equation*}
p(\xi) = \sum_{\substack{\alpha \in \N_0^3, \\ |\alpha| \leq N }} a_\alpha \xi^\alpha
\end{equation*}
be a multi-variate polynomial in $\R^3$ with real coefficients and suppose that $a_\alpha \neq 0$ for some $\alpha \in \mathbb{N}_0^3$ with $|\alpha| = N$. We consider partial differential operators
\begin{equation}
\label{eq:DifferentialOperator}
P(D) = p(-i \nabla_x) = \sum_{|\alpha| \leq N} a_\alpha (-i)^{|\alpha|} \partial^\alpha
\end{equation}
such that for $u \in \mathcal{S}'(\R^3)$ we have
\begin{equation*}
\mathcal{F}( P(D) u) (\xi) = p(\xi) \hat{u}(\xi).
\end{equation*}
By ellipticity we mean that
\begin{equation*}
p_N(\xi) = \sum_{|\alpha| = N} a_\alpha \xi^\alpha \neq 0 
\end{equation*}
for $\xi \neq 0$. We assume $p_N(\xi) > 0$ for the sake of definiteness. In the following we prove existence of solutions $u \in L^q(\R^3)$ such that
\begin{equation*}
P(D) u = f
\end{equation*}
for $f \in L^p(\R^3)$ in a certain range of $p$ and $q$, which satisfy the estimate
\begin{equation*}
\| u \|_{L^q(\R^3)} \lesssim \| f \|_{L^p(\R^3)}.
\end{equation*}
The properties of the vanishing set of $p(\xi)$ play a key role for constructing solutions: Guti\'errez \cite{Gutierrez2004} constructed solutions for $p(\xi) = |\xi|^2 - 1$. In most previous works on elliptic operators was assumed that $\Sigma_0= \{ p(\xi) = 0 \}$ is a smooth manifold with non-vanishing Gaussian curvature $K \neq 0$. In this case the analysis of Guti\'errez applies. Recently, Cast\'eras--F\"oldes \cite{CasterasFoldes2021} analyzed fourth-order Schr\"odinger operators (in dimensions $d \geq 2$) with smooth characteristic surface, and estimates depending on the number of non-vanishing principal curvatures were proved. A wider range was covered in \cite{MandelSchippa2021}, where also surfaces with conic singularities were treated. Presently, we consider the effect of vanishing Gaussian curvature under a transversality assumption, which was described by Erd\H{o}s--Salmhofer \cite{ErdosSalmhofer2007}. The idea of constructing solutions is to consider approximates
\begin{equation*}
\hat{u}_\delta(\xi) = \frac{1}{(2 \pi)^3} \int_{\R^3} \frac{e^{ix.\xi} \hat{f}(\xi)}{p(\xi) + i \delta} d\xi
\end{equation*}
for $\delta \neq 0$ and show uniform bounds
\begin{equation}
\| u_\delta \|_{L^q(\R^3)} \lesssim \| f \|_{L^p(\R^3)}
\end{equation}
for fixed $P(D)$.\\
Then we shall find distributional limits $u \in L^q(\R^3)$, which satisfy
\begin{equation*}
P(D) u = f \text{ in } \mathcal{S}'(\R^3)
\end{equation*}
and
\begin{equation*}
\| u \|_{L^q(\R^3)} \lesssim \| f \|_{L^p(\R^3)}.
\end{equation*}
This is referred to as limiting absorption principle.
We shall still assume that $\nabla p(\xi) \neq 0$ for $\xi \in \Sigma_0$. This is a generic assumption for polynomials. In this case Sokhotsky's formula yields for solutions as described above
\begin{align*}
u(x) &= \frac{1}{(2 \pi)^3} \int_{\R^3} \frac{e^{ix.\xi} \hat{f}(\xi)}{p(\xi) \pm i 0} d\xi \\
 &= \mp \frac{i \pi}{(2 \pi)^3} \int_{\R^3} e^{ix.\xi} \hat{f}(\xi) \delta_{\Sigma_0}(\xi) d\xi + \frac{1}{(2 \pi)^3} v.p. \int_{\R^3} \frac{e^{ix.\xi} \hat{f}(\xi)}{p(\xi)} d\xi.
\end{align*}
This points out a close connection to Fourier restriction. The most basic $L^p$-$L^q$-results rely on the decay of the Fourier transform of the surface measure. This in term is caused by the curvature of the surface. If $K \neq 0$, the estimate
\begin{equation*}
| \hat{\mu}_S(\xi)| = \big| \int_S e^{ix.\xi} dx \big| \lesssim \langle \xi \rangle^{-1}
\end{equation*}
is classical (cf. \cite{Littman1963,Stein1993}). Corresponding $L^p$-$L^q$-estimates for solutions were proved in \cite{MandelSchippa2021}.\\
In this note we consider vanishing total curvature under additional transversality assumptions. For constructing solutions as laid out above, we also have to consider level sets $\Sigma_a = \{ p(\xi) = a \}$ for $|a| \leq \delta_0$. We recall the assumptions in Erd\H{o}s--Salmhofer:\\
Let $I$ be a compact interval and let $\mathcal{D} = e^{-1}(I)$. Suppose that $\Sigma_a$ is a two-dimensional submanifold for each $a \in I$. Let $f \in C^\infty_c(\mathcal{D})$ and define
\begin{equation}
\label{eq:FourierTransformSurface}
\hat{\mu}_a(x) = \int_{\Sigma_a} e^{ix.\xi} f(\xi) d\sigma_a(\xi)
\end{equation}
the Fourier transform of the surface carried measure $f d\sigma_a$.\\
Let $C_0 = \text{diam}(\mathcal{D})$, $C_1 = \| p \|_{C^5(\mathcal{D})}$. The following assumptions have to be met:\\
\textbf{Assumption 1}: 
\begin{equation}
\label{eq:RegularFoliation}
C_2 = \min_{\xi \in \mathcal{D}} |\nabla p(\xi) | > 0,
\end{equation}
which means that $(\Sigma_a)_{a \in I}$ is a regular foliation of $\mathcal{D}$.

Let $K: \mathcal{D} \to \R$ be the Gaussian curvature of the foliation, i.e., for $\xi \in \Sigma_a \subseteq \mathcal{D}$, $K(\xi)$ denotes the curvature of $\Sigma_a$ at $\xi$.\\
The crucial assumption is that the vanishing set of the Gaussian curvature is a submanifold, which intersects $(\Sigma_a)_{a \in I}$ transversally:\\
\textbf{Assumtion 2}: Let $\mathcal{C} = \{ \xi \in \mathcal{D}: K(\xi) = 0 \}$. Then
\begin{equation*}
C_3 = \min_{\xi \in \mathcal{D}} (\{ |\nabla p(\xi) \times \nabla K(\xi)| \, : \xi \in \mathcal{C} \}) >0.
\end{equation*}
With $\nabla K$ non-vanishing on $\mathcal{C}$, it is a two-dimensional submanifold by the regular value theorem. Since $p$ and $K$ are smooth, we find that
\begin{equation*}
\Gamma_a = \mathcal{C} \cap \Sigma_a
\end{equation*}
is a finite union of disjoint regular curves on $\Sigma_a$ for each $a \in I$.\\
Let
\begin{equation*}
\xi \mapsto w(\xi) = \frac{\nabla p(\xi) \times \nabla K(\xi)}{|\nabla p(\xi) \times \nabla K(\xi)|}
\end{equation*}
be the unit vectorfield tangent to $\Gamma_a$. Denote the normal map $\nu: \mathcal{D} \to \mathbb{S}^2$ by
\begin{equation*}
\nu(\xi) = \frac{\nabla p(\xi)}{|\nabla p(\xi)|}.
\end{equation*}
Recall that the Gaussian curvature is given by the Jacobian of the normal map restricted to each surface, $\nu: \Sigma_a \to \mathbb{S}^2$: $K(\xi) = \det \nu'(\xi)$.

We further require the following regularity assumption on the Gauss map.\\
\textbf{Assumption 3}: The number of preimages of $\nu: \Sigma_a \to \mathbb{S}^2$ is finite, i.e.,
\begin{equation*}
C_4 = \sup_{a \in I} \sup_{ \omega \in \mathbb{S}^2} \text{card} \{ p \in \Sigma_a \, : \, \nu(p) =\omega \} < \infty.
\end{equation*}
On the curves $\Gamma_a$, exactly one of the principal curvatures vanish. We define a (local) unit vectorfield $Z \in T \Sigma_a$ along $\Gamma_a$ in the tangent plane of $\Sigma_a$. $Z$ can be extended to a neighbourhood of $\Gamma_a$ as the direction of the principal curvature that is small and vanishes on $\Gamma_a$. We assume that $Z$ is transversal to $\Gamma_a$. Weaker, non-uniform decay estimates were proved in \cite{ErdosSalmhofer2007} also in the presence of tangential points. To ensure uniform decay, we assume the following:\\
\textbf{Assumption 4:} The set of tangential points
\begin{equation*}
\mathcal{T}_a = \{ \xi \in \Gamma_a \, : \, Z(\xi) \times w(\xi) = 0 \},
\end{equation*}
is empty.

\medskip

Under the above assumptions, Erd\H{o}s--Salmhofer \cite[Theorem~2.1]{ErdosSalmhofer2007} proved the following dispersive estimate for the Fourier transform of the surface measure $\mu_a$:
\begin{equation}
\label{eq:DispersiveEstimate}
|\hat{\mu}_a(\xi)| \leq C \langle \xi \rangle^{-\frac{3}{4}}
\end{equation}
with $C=C(C_0,\ldots,C_4,\| f \|_{C^2(\mathcal{D})})$. This morally corresponds to a decay from $\frac{3}{2}$ principal curvatures bounded from below in modulus and thus improves the previous result for one non-vanishing principal curvature (cf. \cite[Theorem~1.3]{MandelSchippa2021}).
In this article we record its consequence for solutions to elliptic differential operators. Allowing for tangential points covers generic surfaces in $\R^3$ as pointed out in \cite{ErdosSalmhofer2007}. However, the decay proved in \cite{ErdosSalmhofer2007} is not uniform in this case anymore. It might be possible to show the same results for a broader class via the estimates due to Ikromov--M\"uller \cite{IkromovMueller2011}.

In the first step, we derive a Fourier restriction--extension theorem for surfaces $\Sigma_a$ by following along the lines of the preceding work \cite{MandelSchippa2021}. We prove strong bounds
\begin{equation}
\label{eq:StrongBoundsFourierRestrictionExtension}
\| \int_{\R^3} e^{ix.\xi} \delta_{\Sigma_a}(\xi) \beta(\xi) \hat{f}(\xi) d\xi \|_{L^q(\R^3)} \lesssim \| f \|_{L^p(\R^3)}
\end{equation}
within a pentagonal region. Here $\beta \in C^\infty_c$ localizes to a suitable neighbourhood of $\{K=0\}$ in $\big( \Sigma_a \big)_{a \in [-\delta_0,\delta_0]}$. Away from $\{ K=0\}$, \cite[Theorem~1.3]{MandelSchippa2021} provides better estimates for $d=3$, $k=2$. On part of the boundary of the pentagonal region, we show weak bounds
\begin{align}
\label{eq:WeakBoundI}
\| \int_{\R^3} e^{ix.\xi} \delta_{\Sigma_a}(\xi) \beta(\xi) \hat{f}(\xi) d\xi \|_{L^{q,\infty}(\R^3)} \lesssim \| f \|_{L^p(\R^3)} \\
\label{eq:WeakBoundII}
\| \int_{\R^3} e^{ix.\xi} \delta_{\Sigma_a}(\xi) \beta(\xi) \hat{f}(\xi) d\xi \|_{L^q(\R^3)} \lesssim \| f \|_{L^{p,1}(\R^3)},
\end{align}
and lastly, restricted weak bounds
\begin{equation}
\label{eq:RestrictedWeakBound}
\| \int_{\R^3} e^{ix.\xi} \delta_{\Sigma_a}(\xi) \beta(\xi) \hat{f}(\xi) d\xi \|_{L^{q,\infty}(\R^3)} \lesssim \| f \|_{L^{p,1}(\R^3)}
\end{equation}
at its inner endpoints. We refer to Figure \ref{fig:FourierRestrictionExtensionEstimate} for a diagram.
For $X,Y \in [0,1]^2$ we write $[X,Y] = \{ Z: \, \exists \lambda \in [0,1]: \, Z = \lambda X + (1-\lambda) Y \}$ and correspondingly $(X,Y)$, $(X,Y]$, etc.

\begin{proposition}
\label{prop:FourierRestrictionExtension}
Let $p:\R^3 \to \R$ be an elliptic polynomial with $\delta_0 > 0$ such that for $\Sigma_a = \{ p(\xi) = a \}$, $-\delta_0 \leq a \leq \delta_0$ Assumptions 1-4 are satisfied in a neighbourhood of $K=0$ in $\Sigma_a$. 
Then, we find \eqref{eq:StrongBoundsFourierRestrictionExtension} to hold for $(\frac{1}{p},\frac{1}{q}) \in [0,1]^2$ provided that
\begin{equation*}
\frac{1}{p} > \frac{7}{10}, \quad \frac{1}{q} < \frac{3}{10}, \quad \frac{1}{p} - \frac{1}{q} \geq \frac{4}{7}.
\end{equation*}

Let
\begin{equation*}
B= \big( \frac{7}{10}, \frac{9}{70} \big), \; C = \big( \frac{7}{10}, 0 \big), \quad B' = \big( \frac{61}{70}, \frac{3}{10} \big), \; C' = \big(1, \frac{3}{10} \big):
\end{equation*}
Furthermore, we find \eqref{eq:WeakBoundI} to hold for $(1/p,1/q) \in (B',C']$, \eqref{eq:WeakBoundII} for $(1/p,1/q) \in (B,C]$, and \eqref{eq:RestrictedWeakBound} for $(1/p,1/q) \in \{B,B' \}$.
\end{proposition}
In the second step we foliate a neighbourhood $U$ of $\Sigma_0$ with level sets of $p$ to show bounds $\| A_\delta f \|_{L^q} \lesssim \| f \|_{L^p(\R^3)}$ for
\begin{equation}
\label{eq:SingularMultiplier}
A_\delta f(x) = \int_{\R^3} \frac{e^{ix.\xi} \beta_1(\xi)}{p(\xi) + i \delta} \hat{f}(\xi) d\xi
\end{equation}
independent of $\delta$. Here, $p,q$ are as in Proposition \ref{prop:FourierRestrictionExtension} and $|p(\xi)| \leq \delta_0$ for $\xi \in \text{supp } (\beta_1)$ with $\Sigma_0 \subseteq \text{supp }(\beta_1)$. Away from the singular set, estimates for
\begin{equation}
\label{eq:ExteriorDomainMultiplier}
B_\delta f(x) = \int_{\R^3} \frac{e^{ix.\xi} \beta_2(\xi)}{p(\xi) + i \delta} \hat{f}(\xi) d\xi
\end{equation}
with $\beta_1 + \beta_2 \equiv 1$ follow from Young's inequality and properties of the Bessel potential. The estimate of $\| B_\delta \|_{L^p \to L^q}$ depends on the order of the elliptic operator.\\
The method of proof is well-known and detailed in \cite{MandelSchippa2021}; see also \cite{KwonLee2020,JeongKwonLee2016} and references therein. We shall be brief. It turns out that one can follow along the lines of \cite{MandelSchippa2021} very closely, substituting $k=\frac{3}{2}$ non-vanishing principal curvatures. We prove the following:
\begin{theorem}
\label{thm:LpLqEstimates}
Let $p: \R^3 \to \R$ be an elliptic polynomial of degree $N \geq 2$. Let $1 < p_1, p_2, q < \infty$ and $f \in L^{p_1}(\R^3) \cap L^{p_2}(\R^3)$. Suppose that there is $\delta_0 > 0$ such that Assumptions 1-4 are satisfied for $( \Sigma_a )_{a \in [-\delta_0,\delta_0]}$. Then, there is $u \in L^q(\R^3)$ satisfying
\begin{equation*}
P(D) u = f
\end{equation*}
in the distributional sense and the estimate
\begin{equation*}
\| u \|_{L^q(\R^3)} \lesssim \| f \|_{L^{p_1}\cap L^{p_2}(\R^3)}
\end{equation*}
holds true provided that
\begin{equation*}
\frac{1}{p_1} > \frac{7}{10}, \quad \frac{1}{q} < \frac{3}{10}, \quad \frac{1}{p_1} - \frac{1}{q} \geq \frac{4}{7}
\end{equation*}
and for $N \leq 3$
\begin{equation*}
0 \leq \frac{1}{p_2} - \frac{1}{q} \leq \frac{N}{3}, \quad \big( \frac{1}{q}, \frac{1}{p_2} \big) \notin \begin{cases}
\{(0,\frac{2}{3}), \; (\frac{1}{3},1)\} \text{ for } N=2, \\
\{(0,1)\} \text{ for } N =3.
\end{cases}
\end{equation*}
\end{theorem}

\section{The Fourier restriction-extension estimate}
\label{section:FourierRestrictionExtensionEstimate}
The purpose of this section is to prove Proposition \ref{prop:FourierRestrictionExtension}. We shall follow the argument of \cite[Section~4]{MandelSchippa2021}. In the first step, we localize to a small neighbourhood of the vanishing set $\{ K = 0 \}$, which by assumptions is a two-dimensional manifold in $\mathcal{D}$. In the complementary set, by compactness, we can apply \cite[Theorem~1.3]{MandelSchippa2021}, which gives uniform $L^p$-$L^q$-estimates in a broader range. Thus, it is enough to suppose that Assumptions 1-4 are valid in a neighbourhood of $\{ K = 0 \}$. The proof follows \cite[Section~4]{MandelSchippa2021} closely.
In the first step, by finite decomposition and rotations, we change to parametric representation of $\Sigma_a = \{(\xi', \psi(\xi')) \, : \xi' \in B(0,c) \}$. We show bounds $T: L^p(\R^3) \to L^q(\R^3)$ for
\begin{equation*}
T f(x) = \int_{\R^3} \delta(\xi_3 - \psi(\xi')) e^{ix.\xi} \chi(\xi') \hat{f}(\xi) d\xi. 
\end{equation*}
The following decay estimate, which is \eqref{eq:DispersiveEstimate}, is central.
\begin{equation*}
\left| \int e^{i(x'.\xi' + x_3 \psi(\xi'))} \beta(\xi') d\xi' \right| \lesssim (1+|x_3|)^{-\frac{3}{4}}.
\end{equation*}
Applying the $TT^*$ argument (cf. \cite{Tomas1975,GinibreVelo1979,KeelTao1998}), we find the following Strichartz estimate:
\begin{equation}
\label{eq:StrichartzEstimate}
\left\| \int e^{i(x'.\xi' + x_3 \psi(\xi'))} \beta(\xi') \hat{f}(\xi') d\xi' \right\|_{L^{\frac{14}{3}}_x(\R^3)} \lesssim \| f \|_{L^2_{\xi'}(B(0,c))}.
\end{equation}

We recall the following lemma to decompose the delta distribution:
\begin{lemma}[{\cite[Lemma~2.1]{ChoKimLeeShim2005}}]
\label{lem:DecompositionLemma}
There is a smooth function $\phi$ satisfying $\text{supp} (\hat{\phi}) \subseteq \{ t \, : \, |t| \sim 1 \}$ such that for all $f \in \mathcal{S}(\R^d)$,
\begin{equation*}
\langle \delta(\xi_3-\psi(\xi')) , f \rangle = \sum_{j \in \Z} 2^{j} \int_{\R^3} \phi(2^j(\xi_3 - \psi(\xi^\prime))) \chi(\xi^\prime) f(\xi) d\xi.
\end{equation*}
\end{lemma}
By this, we can write
\begin{equation*}
T f(x) = \sum_{j \in \Z} 2^j \int_{\R^3} \phi(2^j(\xi_3 - \psi(\xi'))) e^{i x.\xi} \chi(\xi') \hat{f}(\xi) d\xi \\
:= \sum_{j \in \Z} 2^j T_{2^{-j}} f.
\end{equation*}
As pointed out in \cite{ChoKimLeeShim2005}, the contribution of $j \leq 0$ is easier to estimate.

The contribution of $j \geq 0$, i.e., close to the singularity, is estimated by Strichartz and kernel estimates:
\begin{lemma}[{cf. \cite[Lemma~4.3]{MandelSchippa2021}}]
\label{lem:StrichartzEstimate}
Let $q \geq \frac{14}{3}$. Then, we find the following estimate to hold:
\begin{equation*}
\| T_{2^j} f \|_{L^q}(\R^3) \lesssim 2^{\frac{-j}{2}} \| f \|_{L^2(\R^3)}.
\end{equation*}
\end{lemma}
This estimate does not admit summation. For this purpose, we interpolate with the kernel estimate:
\begin{lemma}[{cf. \cite[Lemma~4.4]{MandelSchippa2021}}]
\label{lem:KernelEstimate}
Let
\begin{equation*}
K_\delta(x) = \int_{\R^3} e^{ix.\xi} \beta(\xi') \phi\big( \frac{\xi_3 - \psi(\xi')}{\delta} \big) d\xi.
\end{equation*}
Then $K_\delta$ is supported in $\{(x',x_3):|x_3| \sim \delta^{-1} \}$, and we find the following estimates to hold:
\begin{align*}
|K_\delta(x)| &\lesssim_N \delta^N ( 1 + \delta |x|)^{-N}, \text{ if } |x'| \geq c |x_3|, \\
|K_\delta(x)| &\lesssim \delta^{\frac{7}{4}}, \text{ if } |x'| \leq c|x_3|.
\end{align*}
\end{lemma}
The last ingredient to show (restricted) weak endpoint estimates is Bourgain's summation argument (cf. \cite{Bourgain1985,CarberySeegerWaingerWright1999} and \cite[Lemma~2.3]{Lee2003} for an elementary proof):
\begin{lemma}
\label{lem:SummationLemma}
Let $\varepsilon_1,\varepsilon_2 > 0$, $1 \leq p_1, \, p_2 \leq \infty$, $1 \leq q_1,q_2 < \infty$. For every $j \in \mathbb{Z}$ let $T_j$ be a linear operator, which satisfies
\begin{align*}
\| T_j(f) \|_{q_1} &\leq M_1 2^{\varepsilon_1 j}  \| f \|_{p_1} \\
\| T_j(f) \|_{q_2} &\leq M_2 2^{-\varepsilon_2 j} \| f \|_{p_2}.
\end{align*}
Then, for $\theta, q$ and $p_i$ defined by $\theta = \frac{\varepsilon_2}{\varepsilon_1 \varepsilon_2}$, $\frac{1}{q} = \frac{\theta}{q_1} + \frac{1-\theta}{q_2}$ and $\frac{1}{p} = \frac{\theta}{p_1} + \frac{1-\theta}{p_2}$, the following hold:
\begin{align}
\label{eq:SummationI}
\| \sum_j T_j(f) \|_{q,\infty} &\leq C M_1^\theta M_2^{1-\theta} \| f \|_{p^,1}, \\
\label{eq:SummationII} \| \sum_j T_j(f) \|_q &\leq C M_1^\theta M_2^{1-\theta} \| f \|_{p,1} \text{ if } q_1 = q_2 = q, \\
\label{eq:SummationIII} \| \sum_j T_j(f) \|_{q,\infty} &\leq C M_1^\theta M_2^{1-\theta}  \| f \|_{p} \text{ if } p_1 = p_2.
\end{align}
\end{lemma}

We interpolate the bounds
\begin{equation*}
2^{j} \| T_{2^{-j}} f \|_{L^q(\R^3)} \lesssim 2^{\frac{j}{2}} \| f \|_{L^2(\R^3)}, \quad \frac{14}{3} \leq q \leq \infty,
\end{equation*}
and
\begin{equation*}
2^{j} \| T_{2^{-j}} f \|_{L^\infty(\R^3)} \lesssim 2^{- \frac{3j}{4}} \| f \|_{L^1(\R^3)}
\end{equation*}
as above together with duality to find restricted weak endpoint bounds
\begin{equation*}
\| T f \|_{L^{q,\infty}(\R^3)} \lesssim \| f \|_{L^{p,1}(\R^3)}
\end{equation*}
for $(1/p,1/q) \in \{B, B'\}$, weak bounds
\begin{equation*}
\| T f \|_{L^{q,\infty}} \lesssim \| f \|_{L^p}, \quad \| T f \|_{L^q} \lesssim \| f \|_{L^{p,1}}
\end{equation*}
for $(1/p,1/q) \in (B',C']$, respectively, $(1/p,1/q) \in (B,C]$,
and strong bounds in the interior of the pentagon $\text{conv}(A,B,C,C',B')$ with $A=(1,0)$,
\begin{equation*}
B= \big( \frac{7}{10}, \frac{9}{70} \big), \; C = \big( \frac{7}{10}, 0 \big), \quad B' = \big( \frac{61}{70}, \frac{3}{10} \big), \; C' = \big(1, \frac{3}{10} \big):
\end{equation*}
\begin{center}
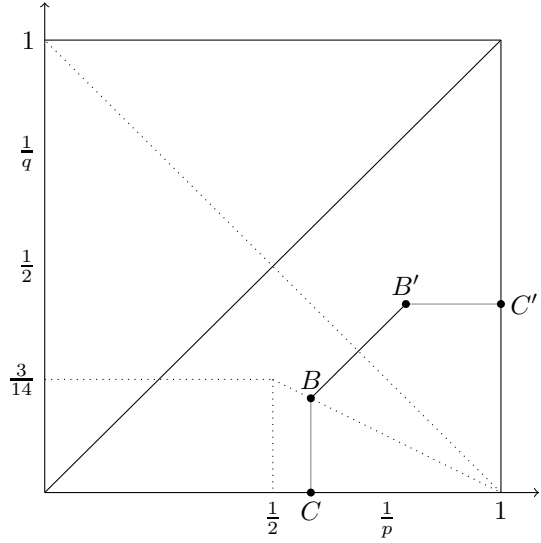
\begin{figure}
\label{fig:FourierRestrictionExtensionEstimate}
\begin{tikzpicture}[scale=0.5]
\draw[->] (0,0) -- (13,0); \draw[->] (0,0) -- (0,13);
\draw (0,0) --(12,12); \draw(0,12) -- (12,12); \draw (12,12) -- (12,0);
\coordinate (E) at (6,3);
\coordinate [label=left:$\frac{3}{14}$] (EX) at (0,3);
\coordinate [label=above:$B$] (B) at (7,2.5);
\coordinate [label=above:$B'$] (B') at (9.5,5);
\coordinate [label=below:$C$] (C) at (7,0);
\coordinate [label=right:$C'$] (C') at (12,5);
\coordinate [label=left:$\frac{1}{2}$] (Y) at (0,6);
\coordinate [label=left:$1$] (YY) at (0,12);
\coordinate [label=below:$\frac{1}{2}$] (X) at (6,0);
\coordinate [label=below:$1$] (XX) at (12,0);
\coordinate [label=left:$\frac{1}{q}$] (YC) at (0,9);
\coordinate [label=below:$\frac{1}{p}$] (XC) at (9,0);

\draw [dotted] (EX) -- (E); \draw [dotted] (E) -- (X); \draw [dotted] (E) -- (XX); \draw [dotted] (XX) -- (YY);
\draw [help lines] (C) -- (B); \draw (B) -- (B');
\draw [help lines] (B') -- (C');

\foreach \point in {(C),(C'),(B),(B')}
	\fill [black, opacity = 1] \point circle (3pt);

\end{tikzpicture}
\caption{Pentagonal region, within which strong $L^p$-$L^q$-Fourier restriction extension estimates hold.}
\end{figure}
\end{center}
Real interpolation of the weak bounds at $B$ and $B'$ gives strong bounds on $(B,B')$. This finishes the proof of Proposition \ref{prop:FourierRestrictionExtension}. \hfill $\Box$

\section{$L^p$-$L^q$-estimates for solutions to elliptic differential operators}

In this section we prove Theorem \ref{thm:LpLqEstimates} relying on Proposition \ref{prop:FourierRestrictionExtension}. The argument parallels \cite[Section~5.2]{MandelSchippa2021} very closely, to avoid repitition we shall be brief. Let $A_\delta$ and $B_\delta$ be as in \eqref{eq:SingularMultiplier} and \eqref{eq:ExteriorDomainMultiplier}. We start with the more difficult estimate of $A_\delta$. We show boundedness of $A_\delta:L^p(\R^3) \to L^q(\R^3)$ independently of $\delta$ with $p$, $q$ as in Proposition \ref{prop:FourierRestrictionExtension}. For this it is enough to show restricted weak type bounds
\begin{equation*}
\| A_\delta \|_{L^{q_0,\infty}} \lesssim \| f \|_{L^{p_0,1}}
\end{equation*}
for $(1/p_0,1/q_0) = (61/70,3/10)$ and the bounds
\begin{equation*}
\| A_\delta f \|_{L^q} \lesssim \| f \|_{L^{p,1}}
\end{equation*}
for $(1/p,1/q) \in ((61/70,3/10),(1,3/10)]$ as strong bounds for $A_\delta$ with $p,q$ as in Proposition \ref{prop:FourierRestrictionExtension} are recovered by interpolation and duality. As $\nabla p(\xi) \neq 0$ for $\xi \in \text{supp} (\beta_1)$ by construction, we can change to generalized polar coordinates. Let $\xi = \xi(p,q)$, where $p$ and $q$ are complementary coordinates.\\
Write
\begin{equation*}
A_\delta f(x) = \int \frac{e^{ix.\xi} \beta_1(\xi)}{p(\xi) + i\delta} \hat{f}(\xi) d\xi = \int dp \int dq \frac{e^{ix.\xi(p,q)} \beta(\xi(p,q)) h(p,q) \hat{f}(\xi(p,q))}{p + i \delta},
\end{equation*}
where $h$ denotes the Jacobian. We can suppose that $|\partial^\alpha h | \lesssim_\alpha 1$ choosing $\text{supp}( \beta)$ small enough. The expression is estimated as in \cite[Subsection~5.2]{MandelSchippa2021} by suitable decompositions in Fourier space and crucially depending on the Fourier restriction estimates for Proposition \ref{prop:FourierRestrictionExtension}; see \cite{KwonLee2020} for $p(\xi) = |\xi|^\alpha$. We write
\begin{equation*}
\frac{1}{p(\xi) + i \delta} = \frac{p(\xi)}{p^2(\xi) + \delta^2} - i \frac{\delta}{p^2(\xi) + \delta^2} = \mathfrak{R}(\xi) - i \mathfrak{I}(\xi).
\end{equation*}
As in \cite{MandelSchippa2021}, $\mathfrak{I}(D)$ is estimated by Minkowski's inequality and Fourier restriction--extension estimates, in the present context from Proposition \ref{prop:FourierRestrictionExtension}. The only difference in the estimate of $\mathfrak{R}(D)$ is that \cite[Lemma~5.1]{MandelSchippa2021} is applied for $k=\frac{3}{2}$ according to the dispersive estimate \eqref{eq:DispersiveEstimate}. For details we refer to \cite[Section~4]{MandelSchippa2021}. This finishes the proof of the estimate for $A_\delta$.

For the estimate of $B_\delta$, we carry out a further decomposition in Fourier space: By ellipticity, there is $R\geq 1$ such that
\begin{equation*}
|p(\xi)| \gtrsim |\xi|^N
\end{equation*}
provided that $|\xi| \geq R$. Let $\beta_2(\xi) = \beta_{21}(\xi) + \beta_{22}(\xi)$ with $\beta_{21}, \beta_{22} \in C^\infty$ and $\beta_{22}(\xi) = 0 $ for $|\xi| \leq R$, $\beta_{22}(\xi) = 1$ for $|\xi| \geq 2R$.\\
We can estimate
\begin{equation*}
\| B_\delta (\beta_{21}(D) f) \|_{L^q} \lesssim \| f \|_{L^p}
\end{equation*}
for any $1 \leq p \leq q \leq \infty$ by Young's inequality uniform in $\delta$. This gives no additional assumptions on $p$ and $q$. We estimate the contribution of $\beta_{22}$ by properties of the Bessel kernel (cf. \cite[Theorem~30]{CossettiMandel2020}) 
\begin{equation*}
\| B_\delta (\beta_{22}(D)  f) \|_{L^q(\R^3)} \lesssim \| \beta_{22}(D) f \|_{L^p(\R^3)}
\end{equation*}
for $1 \leq p,q \leq \infty$ and $0 \leq \frac{1}{p} - \frac{1}{q} \leq \frac{N}{3}$ with the endpoints excluded for $N \leq 3$. For $N \geq 4$ this estimate holds true for $1 \leq p \leq q \leq \infty$. This corresponds to the second assumption on $p$ and $q$ in Theorem \ref{thm:LpLqEstimates}. Lastly, we give the standard argument for constructing solutions: For $\delta > 0$, consider the approximate solutions $u_\delta \in L^q(\R^3)$
\begin{equation*}
\hat{u}_\delta(\xi) = \frac{\hat{f}(\xi)}{p(\xi) + i \delta}.
\end{equation*}
By the above, we have uniform bounds
\begin{equation*}
\| u_\delta \|_{L^q(\R^3)} \lesssim \| f \|_{L^{p_1}(\R^3) \cap L^{p_2}(\R^3)}.
\end{equation*}
By the Banach--Alaoglu--Bourbaki theorem, we find a weak limit $u_\delta \to u$, which satisfies the same bound. We observe that
\begin{equation*}
P(D) u_\delta = f - i \frac{\delta}{P(D)+i \delta} f.
\end{equation*}
Since
\begin{equation*}
\| \frac{\delta}{P(D)+i \delta} f \|_{L^q} \lesssim \delta \| f \|_{L^{p_1} \cap L^{p_2}},
\end{equation*}
we find that $P(D) u_\delta \to f$ in $L^q(\R^3)$. Since $P(D) u_\delta \to P(D) u$ in $\mathcal{S}'(\R^3)$, this shows that
\begin{equation*}
P(D) u = f
\end{equation*}
in $\mathcal{S}'(\R^3)$. The proof is complete. \hfill $\Box$

\section*{Acknowledgements}

Funded by the Deutsche Forschungsgemeinschaft (DFG, German Research Foundation)  Project-ID 258734477 – SFB 1173. I would like to thank Rainer Mandel for discussions on related topics and Jean--Claude Cuenin for pointing out an error in a previous version.

\end{document}